\documentstyle{article}
\input epsf.tex
\begin{document}

\title
{
Center conditions for a simple class of quintic systems
}

\author{
Evgenii P. Volokitin%
\thanks{%
The work is supported by Russian Foundation for Basic Research
}\\
Sobolev Institute of Mathematics, Novosibirsk, 630090, Russia\\
e-mail: volok@math.nsc.ru
}
\date{}
\maketitle

\begin{abstract}
We obtained center conditions for a $O$-symmetric system
of degree 5 for which
the origin is a uniformly isochronous singular point.
\end{abstract}

Classification: Primary 34C05; Secondary 34C25

Keywords: conditions for a center; isochronicity; commutativity

\quad

\noindent {\bf 1.}
Let us consider a planar differential system
\refstepcounter{equation}
\label{1}
$$ \begin{array} {ll}
\dot x = y + x R_{n-1} (x,y),\\
\dot y = - x +y R_{n-1} (x,y),\\
\end{array}
\eqno{(\ref{1})}
$$
where $R_{n-1}(x,y)$ is a polynomial in $x,y$ of degree $n-1.$

System (\ref{1}) has a unique singular point $O(0,0)$ whose linear part of
center type.

Orbits of system (\ref{1}) move around the origin with a constant
angular velocity and the origin is a uniformly isochronous singular point.

In \cite{2} the following problem was proposed:

{\bf Problem 19.1.} {\it Identify systems (\ref{1}) of odd degree
which are $O$-symmetric (not necessarily quasi-homogeneous) having $O$
as a (uniformly isochronous) center.
}

We solve this problem for $n=5$ and derive necessary and sufficient
center conditions for the system
\refstepcounter{equation}
\label{2}
$$ \begin{array} {ll}
\dot x=y+x(a x^2 + b x y + c y^2 + d x^4 + e x^3 y +f x^2 y^2 + g x y^3 +h y^4)%
,\\
\dot y=- x +y(a x^2 + b x y + c y^2 + d x^4 + e x^3 y +f x^2 y^2 + g x y^3 +h y^4)%
,\\
a, b, c, d, e, f, g, h \in R.
\end{array}
\eqno{(\ref{2})}
$$

{\bf Theorem.} {\it The origin is a center of system (\ref{2}) if and
only if one of the following sets of conditions is satisfied:}
\refstepcounter{equation}
\label{3}
$$
\begin{array}{lll}
{\rm (i)}&  a=b=c=0, f=-3(d+h); \hfill \\
{\rm (ii)}&  a=c=d=f=h=0;\\
{\rm (iii)}& a\neq 0, c=-a, f=3b(ae-bd)/(2a^2),\\
 &g=(2a^2bd+(2a^2-b^2)(bd-ae))/(2a^3),\\
 &h=(-2a^2d+b(bd-ae))/(2a^2).
\end{array}
\eqno{(\ref{3})}
$$

{\bf Proof.}

Necessity:

To describe the behaviour of trajectories of (\ref{2}) near the origin
we construct the comparison function \cite{3}
$$
F(x,y)=(x^2+y^2)/2 + f_{3}(x,y) + f_{4}(x,y) + \ldots ,
$$
where $f_k$ is a homogeneous polynomial of degree $k$
whose derivative is
$$ \begin{array}{ll}
\frac{dF}{dt} = D_1(x^4+y^4) + D_2 (x^6 + y^6) + D_3 (x^8 +y^8)  + \ldots
\end{array}
$$

The number of the first coefficient $D_i$ other than zero
defines the multiplicity of a complex focus and the sign of this
coefficient defines stability of a focus; if $D_i=0$ for all $i$
the origin is a center of (\ref{2}). We refer to coefficients $D_i$
as the Poincar\'{e}--Lyapunov constants.

To find the Poincar\'{e}--Lyapunov constants of a system
$\dot x=p(x,y), \dot y=q(x,y)$ with a linear center
we used computer algebra and wrote a {\it Mathematica} code
that rests on the Poincar\'{e} algorithm in \cite{3}; see \cite{4}
for more details.

\begin{verbatim}
PLconst[n_] :=
Module[{dF, ff, fF, x, y, pP, qQ, dD},
       fF[2] := (x^2+y^2)/2;
       fF[i_] := Sum[ff[i-j, j]*x^(i-j)*y^j, {j, 0, i}];
       pP[1] := y;
       pP[i_] := Sum[p[i-j, j]*x^(i-j)*y^j, {j, 0, i}];
       qQ[1] := -x;
       qQ[i_] := Sum[q[i-j, j]*x^(i-j)*y^j, {j, 0, i}];
       dF[k_] := (Sum[D[fF[i], x]*pP[k+1-i], {i, 2, k}]+
          Sum[D[fF[i], y]*qQ[k+1-i], {i, 2, k}])//Expand;
       Do[
          Solve[Table[Coefficient[dF[k], x^(k-j) y^j],
                {j, 0, k}]
                ==Table[0, {k+1}],
                Table[ff[k-j, j], {j, 0, k}]
               ]/.Rule->Set;
          Solve[Table[Coefficient[dF[k+1], x^(k+1-j)*y^j],
                {j, 0, k+1}]
                ==Flatten[{dD[k], Table[0, {k}], dD[k]}]
                &&ff[0, k+1]==0,
                Flatten[{Table[ff[k+1-j, j],{j, 0, k+1}], dD[k]}]
               ]/.Rule->Set,
          {k, 3, 2n+1, 2}];
          Table[Numerator[Together[dD[k]]],{k, 3, 2n+1, 2}]
      ]
\end{verbatim}

The procedure {\tt PLconst[n]} returns a list $\{D_1,\ldots,D_n\}$ of the
Poincar\'{e}-Lyapunov constants if we define the coefficients $p_{ij}, q_{ij}$
$(2\leq i+j \leq 2n+1)$  in the Taylor series expansion
of functions $p(x,y), q(x,y)$ beforehand.

Using this procedure, we found the first four Poincar\'{e}-Lyapunov
constants of (\ref{2}).
$$
\begin{array}{ll}
D_1=&2(a + c),\\
D_2=&-4ab - 4bc + 3d + f + 3h,\\
D_3=&2(-85a^3 + 15ab^2 - 67a^2c + 15b^2c + 61ac^2 + 43c^3-24bd - 34ae -   \\
    & -22ce- 12bf - 50ag - 38cg - 48bh),\\
D_4=&44600a^3b + 2736ab^3 + 84696a^2bc + 2736b^3c + 47688abc^2 +7592bc^3- \\
    & - 37120a^2d- 1782b^2d - 32552acd - 2704c^2d + 2364abe + 1284bce- \\
    & - 2673de - 6120a^2f- 234b^2f -3384acf + 792c^2f - 891ef +  \\
    &+ 6876abg+5076bcg - 3807dg -1269fg + 4720a^2h + 1098b^2h +\\
    & + 31448ach + 19456c^2h -2673eh - 3807gh.
\end{array}
$$

It is easy to verify that
the equalities $D_i=0, i=1,2,3,4$
are equivalent to the following relations
$$
\begin{array}{ll}
a + c=0,&
3d + f + 3h=0,\\
3ce - bf + 3cg - 6bh=0,&
2c^2f - 3bcg + 3b^2h=0.
\end{array}
$$

If $a=0$ then our simultaneous polynomial equations have two sets of
solutions indicated in (i) and (ii). If $a\neq0$
then, in view of the condition $c=-a$, we see that
the other three equations constitute a non degenerate linear
system for determining the variables $f,g,h$. The solution is given by (iii).

The necessity part of the theorem is proved.

Sufficiency:

Case (i). System (\ref{2}) now takes the form
\refstepcounter{equation}
\label{4}
$$
\begin{array}{ll}
\dot x=y+x(d x^4 + e x^3 y +f x^2 y^2 + g x y^3 +h y^4)\equiv y+xp_4(x,y),\\
\dot y=- x +y(d x^4 + e x^3 y +f x^2 y^2 + g x y^3 +h y^4)\equiv x+yp_4(x,y).\\
\end{array}
\eqno{(\ref{4})}
$$

This is a quasi-homogeneous system of degree 5 whose coefficients
satisfy the equality $f=-3(d+h)$ which is the necessary and sufficient
center condition in the case we study \cite{5}.

Case (ii). System (\ref{2}) now takes the form
\refstepcounter{equation}
\label{5}
$$
\begin{array}{ll}
\dot x=y+x^2 y (b + e x^2 + g y^2),\\
\dot y=-x+x y^2 (b + e x^2 + g y^2).
\end{array}
\eqno{(\ref{5})}
$$

The planar differential system
\refstepcounter{equation}
\label{6}
$$
\dot x=p(x,y), \dot y=q(x,y)
\eqno{(\ref{6})}
$$
is said to be reversible (in the sense of \.{Z}ol\c{a}dek) if its orbits are
symmetric with respect to a line passing through the origin.

System (\ref{6}) is reversible if there is a linear transformation
$S:R^2 \to R^2,$ sending a point $(x,y)$ to the point $(x',y')$ symmetric
to $(x,y)$ with respect to the line $\alpha x + \beta y =0$ and
satisfying the condition
$S (p(x,y),q(x,y))=-(p(S(x,y)), q(S(x,y)))$.

A more general condition of reversibility is as follows
$$
\begin{array}{cc}
2 \alpha \beta (p(x,y)p(x',y')-q(x,y)q(x',y'))+\\
(\beta^2-\alpha^2)(p(x,y)q(x',y')+p(x',y')q(x,y))=0.
\end{array}
$$

It is well known that if system (\ref{6}) is reversible and has a linear
center  at the origin then the origin is a center of this system
(see \cite{3}, for example).

Obviously, system (\ref{5}) is reversible because its trajectories
are symmetric with respect to both coordinate axes. So, the origin
is a center for system (\ref{5}).

Case (iii). System (\ref{2}) now takes the form
\refstepcounter{equation}
\label{7}
$$
\begin{array}{ll}
(2a^3) \dot x=&(2a^3) y + x(a x^2  + b x y - a y^2)\\
 &(2 a^3  + 2 a^2 d x^2 - 2 a b d x y + 2 a^2 e x y + 2 a^2 d y^2  - b^2 d y^2 +a b e y^2),\\
(2a^3) \dot y=& -(2a^3)x + y(a x^2  + b x y - a y^2)\\
 &(2 a^3  + 2 a^2 d x^2 - 2 a b d x y + 2 a^2 e x y + 2 a^2 d y^2  - b^2 d y^2 +a b e y^2)\\
\end{array}
\eqno{(\ref{7})}
$$

It turns out that system (\ref{7}) is reversible. Its trajectories
are symmetric with respect to each of the two perpendicular lines
defined by the equation $a x^2 + b x y - a y^2= 0.$
The appropriate linear transformation $S$ is given by each of the
two matrices
$$
S_{1,2}=\pm(4a^2+b^2)^{-1/2}
\left(
\begin{array}{cc}
      -b & 2a \\
      2a & b \\
\end{array}
\right).
$$

This fact is confirmed by the straight calculations. We used
{\it Mathematica} here.

With the coordinate change
$x\mapsto x \cos \varphi + y \sin \varphi,
 y\mapsto -x \sin \varphi + y \cos \varphi$
where the angle $\varphi$ is defined from the condition
$a \tan^2\varphi +b \tan\varphi-a=0,$
system (\ref{7}) becomes as follows
$$
\begin{array}{ll}
\dot x=y+x^2 y (b_1 + e_1 x^2 + g_1 y^2),\\
\dot y=-x+x y^2 (b_1 + e_1 x^2 + g_1 y^2).
\end{array}
$$

Hence the origin is a center for system (\ref{2}) in this case once again.

The theorem is proved.

\quad

\noindent {\bf 2.}
It is known that isochronism of a center of a planar polynomial
system is equivalent to existence of an analytic transversal system
commuting with a given system in a neighbourhood of a center \cite{6};
observe that an arbitrary polynomial system with
isochronous center not necessarily commutes with a polynomial
system \cite{7,8}.

It is proved in \cite{9} that
if the systems
\refstepcounter{equation}
\label{8}
$$
\begin{array}{ll}
\dot x=p(x,y), \dot y=q(x,y)\\
\dot x=r(x,y), \dot y=s(x,y)
\end{array}
\eqno{(\ref{8})}
$$
commute then $\mu(x,y)=1/(p(x,y)s(x,y)-q(x,y)r(x,y))$
is an integrating factor of both systems.

Thereby if both commuting systems are polynomial then we can
find the integrating Darboux factor for the given system and
integrate the latter
(about the method of Darboux and the relevant definitions see \cite{10},
for example).

We now state the following fact which will be useful later.

Considering (\ref{8}), assume that
$$
\begin{array}{ll}
p(x,y)=y+xR(x,y), q(x,y)=-x+yR(x,y),\\
r(x,y)=xQ(x,y), s(x,y)=yQ(x,y),
\end{array}
$$
where $R(x,y), Q(x,y)$ are polynomials in $x,y.$
Then the algebraic curves $x^2+y^2=0, Q(x,y)=0$ are invariants
for each of these systems.

Indeed, it is immediately obvious that $x^2+y^2=0$ is an invariant of
both systems with the cofactor $2R(x,y)$ and $2Q(x,y)$ respectively.
The curve $Q(x,y)=0$ is an invariant of the second system with the cofactor
$x Q_x(x,y)+y Q_y(x,y)$.

Because our systems commute
the Lie bracket of the vector fields $(p, q)$  and $(r, s)$ vanishes and
we have
$$
p_x(x,y) r(x,y) + p_y(x,y) s(x,y)- r_x(x,y) p(x,y) - r_y(x,y) q(x,y) =0,
$$
or
$$
\begin{array}{l}
x Q(x,y)(R(x,y)+xR_x(x,y))+yQ(x,y)(1+xR_y(x,y))-\\
-p(x,y)(Q(x,y)+xQ_x(x,y))-xq(x,y)Q_y(x,y)=0,
\end{array}
$$
or
$$
\begin{array}{l}
x(Q_x(x,y)p(x,y)+Q_y(x,y)q(x,y))=\\
 =(R(x,y)+xR_x(x,y))x Q(x,y)+\\
 +(1+xR_y(x,y))yQ(x,y)-Q(x,y)p(x,y)= \\
 =(R(x,y)+xR_x(x,y))x Q(x,y)+ \\
 +(1+xR_y(x,y))yQ(x,y)-Q(x,y)(y+xR(x,y))= \\
\phantom{x(Q_x(x,y))} =x(xR_x(x,y)+yR_y(x,y))Q(x,y).
\end{array}
$$

We see that the curve $Q(x,y)=0$ is an invariant with the cofactor
$xR_x(x,y)+yR_y(x,y)$.

In this case $\mu(x,y)=1/(Q(x,y)(x^2+y^2))$
is an integrating Darboux factor.

\quad

\noindent {\bf 3.}
In each of the three cases we have found a non trivial polynomial
system commuting with the respective system.

In case (i) such a system is
$$
\begin{array}{ll}
\dot x=x(1+e x^4 -4d x^3 y +4h x y^3 -g y^4)\equiv x(1+ q_4(x,y)),\\
\dot y=y(1+e x^4 -4d x^3 y +4h x y^3 -g y^4)\equiv y(1+ q_4(x,y)).
\end{array}
$$

The function
$$
\mu(x,y)=\frac{1}{(x^2+y^2)(1+q_4(x,y))},
$$
is the integrating Darboux factor of (\ref{4})
and the function
\refstepcounter{equation}
\label{9}
$$
H(x,y)=\frac{(x^2+y^2)^2}{1+q_4(x,y)}
\eqno{(\ref{9})}
$$
is the first rational
integral of (\ref{4}).

The algebraic curves $x^2+y^2=0, 1+q_4(x,y)=0$ are invariant curves
for (\ref{4}).

According to \cite{5} system (\ref{1}) has a center of type
$B^k, 1\leq k \leq n-1$ whose boundary is a finite union of
$k$
unbounded open trajectories. Using (\ref{9}), in case (i) we can
describe this boundary explicitly:
$$
\varrho=\frac{1}{(á_0-q_4(\cos \varphi, \sin \varphi))^{1/4}},
$$
where $c_0=\max_{[0,2\pi]} q_4(\cos \varphi, \sin \varphi),
x=\varrho \cos \varphi, y=\varrho \sin \varphi.$

A straight analysis of this expression allows us to conclude that
in our case a center may be of type $B^2$ or $B^4$ only.

In case (ii) system (\ref{5}) commutes with the system
$$
\begin{array}{ll}
\dot x=(e-g)x +x(ex^2+gy^2)(b+ex^2+gy^2),\\
\dot y=(e-g)y +y(ex^2+gy^2)(b+ex^2+gy^2).
\end{array}
$$

This permits us to find an integrating Darboux factor
\refstepcounter{equation}
\label{10}
$$
\mu(x,y)=\frac{1}{(x^2+y^2)(e-g+(ex^2+gy^2)(b+ex^2+gy^2))}.
\eqno{(\ref{10})}
$$

The algebraic curves $x^2+y^2=0, e-g+(ex^2+gy^2)(b+ex^2+gy^2)=0$
are invariant ones for system (\ref{5}).

If $b=0$ system (\ref{5}) is a system of the form (\ref{4})
for which the condition $f=-3(d+h)$ is obviously fulfilled.
Then its first integral is
$$
H(x,y)=\frac{(x^2+y^2)^2}{1+ex^4-gy^4}.
$$

If $b\neq0$ then we may suppose that $b=1.$ The general case reduces to
this by the change of variables
$x \to x/\sqrt{b}, y \to y/\sqrt{b}$ for $b>0$ or
$x \to y/\sqrt{-b}, y \to x/\sqrt{-b}, t \to -t$ for $b<0.$

%\exp[{-}2{\arctan\frac{1{+}2ex^2{+}2gy^2}{\sqrt{4(e{-}g)-1}}
%/\sqrt{4(e{-}g){-}1}].\\

Then our system takes the form
\refstepcounter{equation}
\label{11}
$$
\begin{array}{ll}
\dot x=y+x^2 y (1 + e x^2 + g y^2)\equiv X_1(x,y),\\
\dot y=-x+x y^2 (1 + e x^2 + g y^2)\equiv Y_1(x,y).
\end{array}
\eqno{(\ref{11})}
$$

The functions
$$
C_1=x^2+y^2,
C_2=e-g + (e x^2 +g y^2) +(e x^2 +g y^2)^2
$$
are invariants for (\ref{11}) with the cofactors
$$
L_1=2xy(1+e x^2+ g y^2),
L_2=2xy(1+2(e x^2+g y^2)).
$$
Moreover, if $e\neq g$ the function
$$
C_3=\exp(\int_0^{e x^2+g y^2} \frac{dt}{e-g+t+t^2})
$$
is invariant with the cofactor $L_3=2xy$.

We have $2L_1-L_2-L_3=0.$ Then the function
$$
H(x,y)=\frac{C_1^2}{C_2 C_3}=\frac{(x^2+y^2)^2}
{(e-g + (e x^2 + g y^2) +(e x^2 + g y^2)^2)
\exp(\int_0^{e x^2+g y^2} \frac{dt}{e-g+t+t^2})}
$$
is the first Darboux integral of (\ref{11}).

Let us remark that
$$
\int \frac{dt}{e-g+t+t^2}=
\frac{2}{\sqrt{4(e-g)-1}}\arctan\frac{1+2t}{\sqrt{4(e-g)-1}}
$$
for $4(e-g)-1 >0$ and
$$
\int \frac{dt}{e-g+t+t^2}=
-\frac{2}{1+2t}
$$
for $4(e-g)-1 =0$ and
$$
\int \frac{dt}{e-g+t+t^2}=
\frac{1}{\sqrt{1-4(e-g)}}\ln\frac{1+2t-\sqrt{1-4(e-g)}}{1+2t+\sqrt{1-4(e-g)}}
$$
for $4(e-g)-1 <0$.

If $e=g$ the function
$$
C_3=\exp^{\frac{1+x^2}{x^2+y^2}}
$$
is invariant with the cofactor for (\ref{11}) $L_3=-2exy$.

We have $2L_1-L_2+\frac{1}{e}L_3=0.$ Then the function
$$
%\begin{array}{l}
H(x,y)=\frac{C_1^2}{C_2 C_3^{1/e}}=\frac{(x^2+y^2)}
{(1+ e x^2 + e y^2)}
\exp(\frac{1}{e}\frac{1+x^2}{x^2+y^2})
%\end{array}
$$
is the first Darboux integral of (\ref{11}) for $e=g$.

Since (\ref{11}) has a unique finite singular point at the origin,
the phase portraits are obtained by studying the points at infinity.
A standard inspection
of the location and types of such points on the equator
of the Poincar\'{e} sphere allows us to conclude that (\ref{11})
has phase portraits of two types only: a center is of type
$B^2$ when $eg\geq0$ or of type $B^4$ when $eg<0.$

In case (iii) a commuting system and first integral may be found on
considering that system (\ref{7}) is equivalent to (\ref{11}).

Observe that for $d=e=0$ system (\ref{7}) is a quasi-homogeneous
$O$-symmetric cubic system of the form
$$
\begin{array}{ll}
\dot x=y+x(ax^2+bxy-ay^2),\\
\dot y=-x+y(ax^2 +bxy-ay^2).
\end{array}
$$
It commutes with the system
$$
\begin{array}{ll}
\dot x=x+x(bx^2-2axy),\\
\dot y=y+y(bx^2-2axy),
\end{array}
$$
and has the first integral
$$
H(x,y)=\frac{x^2+y^2}{1+bx^2-2axy}.
$$

Summarizing we conclude that
the system under consideration has phase portraits of two types only:
a center is of type $B^2$ or of type $B^4$.

\quad


\begin{thebibliography}{10}


\bibitem{2}
{\it R.~Conti,}
Centers of planar polynomial systems. A review.
Le Mathematiche. 1998. Vol.LIII. Fasc.II. P.207--240.

\bibitem{3}
{\it V.~V.~Nemytskii} and {\it V.~V.~Stepanov},
Qualitative Theory of Differential Equations.
Princeton Univ. Press, 1960.

\bibitem{4}
{\it E.~P.~Volokitin} and {\it S.~A.~Treskov},
About the Lyapunov values of a complex focus of a planar dynamical system.
Izv. Ross. Akad. Estestv. Nauk, Mat. Mat. Model. Inform. Upr.
1997. V.11. No. 1. P. 59--72 (in Russian).

\bibitem{5}
{\it R.~Conti,}
Uniformly isochronous centers of polynomial systems in $R^2$.
Differential Equations, Dynamical Systems, and Control Science
(K. D. Elworthy et al., eds),
Lecture Notes in Pure and Appl. Math., vol. 152,
Marcel Dekker, New York, 1994, pp. 21--31.

\bibitem{6}
{\it M.~Sabatini,}
Characterizing isochronous centres by Lie brackets.
Differential Equations and Dynamical Systems. 1997. V. 5. No. 1. P.~91--99.

\bibitem{7}
{\it J.~Devlin,}
Coexisting isochronous and nonisochronous centres,
Bull. London Math. Soc. 1996. V. 28. No.~134. P. 495--500.

\bibitem{8}
{\it F.~P.~Volokitin} and {\it V.~V.~Ivanov},
Isochronicity and commutation of polynomial vector fields.
Siberian Mathematical Journal. 1999. V. 40. No. 1. P. 23--38.

\bibitem{9}
{\it J.~Chavarriga, H.~Giacomini} and {\it J.~Gin\'{e}},
The null divergence factor.
Publ. Mat. 1997. V. 41. P. 41--56.

\bibitem{10}
{\it P.~Marde\v{s}i\'{c}, C.~Rousseau} and {\it B.~Toni},
Linearization of isochronous centers.
J. Differential Equations. 1995. V. 121. No. 1. P. 67--108.

\end{thebibliography}
\end{document}